\documentclass[10pt]{amsart}
\usepackage[all]{xy}
\usepackage{amsmath}
\usepackage{amssymb}
\usepackage{amsthm}
\input hart.sty         
\renewcommand{\>}{\rightarrow}

\begin{document}
\theoremstyle{plain}
\newtheorem{thm}{Theorem}[section]
\newtheorem{prop}[thm]{Proposition}
\newtheorem{cor}[thm]{Corollary}
\newtheorem{lem}[thm]{Lemma}

\theoremstyle{definition}
\newtheorem{defn}[thm]{Definition}
\newtheorem{que}[thm]{Question}
\newtheorem{cla}[thm]{Claim}
\newtheorem{exam}[thm]{Example}
\newtheorem{obs}[thm]{Observation}
\newtheorem{rmk}[thm]{Remark}

\theoremstyle{remark}
\newtheorem{note}[thm]{Note}

\title{Some Inequalities for Kodaira-Iitaka Dimension on Subvarieties} 
\author{Travis Kopp}    
\date{\today}        
\maketitle
\setcounter{section}{-1}
\section{Introduction}
In this paper we will be principally interested in looking at relations between 
the Kodaira-Iitaka dimension of a divisor on an ambient normal variety and the Kodaira-Iitaka dimension 
of related divisors on irreducible normal subvarieties.\\
\indent
We will begin with
some preliminary discussion about reflexive sheaves and Kodaira-Iitaka dimension
in section \ref{Prelim}.
Then in section \ref{Part I}
we will be primarily interested in the question
of how the Kodaira-Iitaka dimension of a divisor on a normal variety compares to
the Kodaira-Iitaka dimension of related divisors on an irreducible normal subvariety of codimension $1$.  
In particular, we will use geometric techniques to show the following,
\begin{thm}\label{intro main}
Let $X$ be a complete normal variety, $Y\sub X$ a normal closed subvariety 
of codimension $1$ and $\sL$ an invertible sheaf on $X$.
Then there exist integers $n_1>0, n_2\geq 0$ such that,
\[
\ka(X,\sL)-1\leq\ka(Y, \sL^{n_1}(-n_2 Y)\vert_Y) \]
Furthermore if $Y$ is not contained in the stable base locus of $\sL$, we may take $n_1\gg n_2$.
\end{thm}
As an application we will show that Theorem \ref{intro main} implies the following inequality
when we are concerned with a subvariety of arbitrary codimension rather than a divisor only.
\begin{thm}\label{intro subvar}
Let $X$ be a complete smooth variety, $A\sub X$ a smooth closed subvariety
and $\sL$ an invertible sheaf on $X$.
Then there exists a constant $\al>0$ and integers $n_1>0,\ n_2\geq0$,
such that the following inequality holds for $t$ sufficiently large and divisible,
\[
\al t^{\ka(X,\sL)-1}< h^0(A,\sL^{tn_1}\vert_A\otimes\sym^{tn_2}\sN_{A\vert X}^*) \]
Furthermore if $A$ is not contained in the stable base locus of $\sL$, we may choose
$n_1\gg n_2$.
\end{thm}

In section \ref{Part II} we will consider the situation of two subvarieties,
$A$ and $B$, in our ambient variety $X$ with $\dim A+\dim B\geq \dim X$.
We will ask that each of these subvarieties moves as a cycle
in a family covering $X$ and that these families have a ``nice'' intersection property.
In particular, we will prove the following.

\begin{thm}\label{intro main2}
Let $X$ be a projective normal variety 
in characteristic $0$
and let $A, B\sub X$  be normal closed subvarieties.
Suppose that the cycles $[A]$ and $[B]$ 
each move in a well defined integral family of proper algebraic cycles
of $X$ over a variety
in the sense of \cite[I.3.10]{Koll} 
and that for a general closed point $x\in X$
there are members of these families, $[A']$ and $[B']$,
such that $A'$ and $B'$ intersect properly and $x\in A'\cap B'$.
Then for any invertible sheaf $\sL$ on $X$,
\[
\ka(X,\sL)\leq \ka(A,\sL\vert_A)+\ka(B,\sL\vert_B) \]
\end{thm} 

The origin of these kinds of questions may be taken to be 
the following well known theorem in the case of a smooth fiber of a proper fibration.

\begin{thm}[Easy Addition, see e.g. \protect{\cite[\S 10]{Iitaka}}]
Let $f:X\>Y$ be a dominant proper morphism of smooth varieties
with connected fibers in characteristic $0$
and let $\sL$ be an invertible sheaf on $X$.
If $Y_{\sm}\sub Y$ is an open set over which $f$ is smooth, then,
\[
\ka(X,\sL)\leq\ka(X_y, \sL_y)+\dim Y\quad \text{for all}\ y\in Y_{\sm}, \]
\end{thm}
In \cite{PSS}
T. Peternell, M. Schneider and A.J. Sommese showed that this could be generalized
to the case of a subvariety, $A\sub X$,
not necessarily a fiber of a fibration, but with some positivity conditions on the normal bundle,
$\sN_{A\vert X}$.\\ 
\indent
For example, they proved the following inequality of Kodaira-Iitaka dimensions,

\begin{thm}[\protect{\cite[Theorem 4.1]{PSS}}] \label{PSS mt}
Let $X$ be a complete smooth complex variety, $A\sub X$ a smooth closed subvariety
and $\sL$ an invertible sheaf on $X$.
If $\sN_{A\vert X}$ is $\Q$-effective
(i.e. $\sym^{m}\sN_{A\vert X}$ is generically generated by global sections for some $m>0$),
then,
\[
\ka(X,\sL)\leq \ka(A,\sL\vert_A)+\codim_X A \]
\end{thm}
This result and other similar ones in \cite{PSS} follow from Theorem \ref{intro subvar}.\\
\indent
The condition in theorem \ref{PSS mt} that $\sN_{A\vert X}$ be $\Q$-effective
can be shown to be met when $A$ moves as a cycle in a covering family of $X$
\cite[Theorem 4.16]{PSS}.
Thus theorem \ref{intro main2} may also be seen as a natural generalization of this result 
to the case of a pair of subvarieties.\\
\indent
The results in this paper extend those in \cite{Kopp1} and are also presented in my doctoral thesis at the University of Washington.
I would like to thank my advisor, S\'{a}ndor Kov\'{a}cs, for his very valuable advice and 
encouragement throughout the process of preparing and writing this paper.

\section{Preliminary Discussion of Kodaira-Iitaka Dimension}\label{Prelim}
\subsection{Reflexive Sheaves, Rational Maps, and Kodaira-Iitaka Dimension}\ \\
\label{RS}
\indent
Throughout this paper ``variety'' will be taken to mean an integral separated scheme of finite type
over an algebraically closed field. 
In section \ref{Part II} we will also ask that the field be of characteristic $0$.\\
\indent
We will begin with some preliminary observations
concerning reflexive sheaves, rational maps, and Kodaira-Iitaka dimension.
Using reflexive sheaves will allow the statement of some of our results to be more general
than would be possible if we restricted ourselves to invertible sheaves only.\\
\indent
A coherent sheaf $\sF$ on a scheme $X$ is called \textit{reflexive} if
the natural map, $\sF\>\sF^{**}$, to the double dual is an isomorphism.
On an integral noetherian scheme the dual of every coherent sheaf is reflexive.
We will be interested in the reflexive symmetric power and reflexive tensor product,
defined as follows,
\begin{defn}
For an integral noetherian scheme $X$ 
and reflexive sheaves $\sF$ and $\sG$ on $X$, let,
\begin{enumerate}
\item
$\sym^{[\, t]}\sF :=(\sym^t\sF)^{**}$, written as $\sL^{[\, t]}$
for a rank one reflexive sheaf $\sL$.
\item 
$\sF[\otimes]\sG :=(\sF\otimes\sG)^{**}$
\end{enumerate}
\end{defn}
We will also want to use the following properties of reflexive sheaves.
Many of these are found in \cite[\S 1]{Hart2}.
\begin{prop}\ 
\begin{enumerate}
\item
The dual of any coherent sheaf on an integral noetherian scheme is reflexive.
\item
If $X$ is a normal integral noetherian scheme, $\sF$ a reflexive sheaf on $X$, and $U\sub X$ an open set with
$\codim_X(X\setminus U)\geq 2$, then $\sF=i_*(\sF\vert_U)$, where $i:U\>X$ is the inclusion.
Thus $H^0(X,\sF)\iso H^0(U,\sF\vert_U)$.
\item
If $X$ is a normal noetherian scheme, $U$ as before, and $\sF$ a reflexive sheaf on $U$, then
$i_*\sF$ is a reflexive sheaf on $X$ of the same rank.
\item
If $X$ is a regular noetherian scheme, then any rank one reflexive sheaf on $X$ is an invertible sheaf.
\end{enumerate}
\end{prop}
We will also want some definitions to define a general class of varieties with which we'll work.
The following is inspired by and compatible with the definition of a small normal pair found in \cite{PSS}.
\begin{defn}\ 
\begin{enumerate}
\item
A normal variety $X$ will be called \textit{almost complete} if 
there is an open immersion $X\sub\overline{X}$ of $X$ 
as an open set in a complete normal variety $\overline{X}$ 
with $\codim_{\overline{X}}(\overline{X}\setminus X)\geq 2$.
\item
A pair $(X,A)$ will be called a \textit{small normal pair}, if both $A$ and $X$ are almost complete normal varieties,
$A\sub X$, and $\codim_A(A\cap\Sing X)\geq 2$.
\end{enumerate}
\end{defn}
Let $X$ be a normal variety and $X_{\sm}\sub X$ be its smooth locus.
Then we have $\codim_X(X\setminus X_{\sm})\geq 2$.
If $\sL$ is a rank one reflexive sheaf on $X$ then $\sL\vert_{X_{\sm}}$ is an invertible sheaf on $X_{\sm}$
with $H^0(X,\sL)\iso H^0(X_{\sm},\sL\vert_{X_{\sm}})$.\\
\indent
This means that given a non-trivial finite dimensional linear subspace $L\sub H^0(X,\sL)$,
we may naturally define a subspace of $H^0(X_{\sm},\sL\vert_{X_{\sm}})$,
then a rational map $X_{\sm}\ratmap\P^N$ up to linear isomorphism of $\P^N$,
and finally a rational map $\phz_L:X\ratmap\P^N$, also up to linear isomorphism of $\P^N$.
In the case that $L=H^0(X,\sL)$ we will write $\phz_{\sL}$ for $\phz_L$.\\
\indent
Conversely, suppose $\phz:X\ratmap\P^N$ is a rational map
such that $\phz(X)$ is not contained in any proper linear subspace of $\P^N$.
Let $U$ be the intersection of the smooth locus of $X$ and the locus where $\phz$ is defined.
Then $\codim_X(X\setminus U)\geq 2$ and
$\phz\vert_U$ is determined by a linear subspace of $H^0(U,\phz\vert_U^*\sO_{\P^N}(1))$.
This means $\phz=\phz_L$ for some $L\sub H^0(X,\phz^*\sO_{\P^N}(1))$,
where $\phz^*\sO_{\P^N}(1)=i_*(\phz\vert_U^*\sO_{\P^N}(1))$ is a rank one reflexive sheaf on $X$.
In general if $X\ratmap V$ is a rational map defined on $U\sub X$ and $\sL$ is an invertible sheaf on $V$,
we will write $\phz^*\sL$ for the rank one reflexive sheaf $i_*(\phz\vert_U^*\sL)$.  
\\
\indent
This allows us to define a notion of Kodaira-Iitaka dimension for rank one reflexive sheaves on almost complete normal varieties,
\begin{defn}[Kodaira-Iitaka Dimension]\label{KID def}
Let $X$ be an almost complete normal variety and $\sL$ a rank one reflexive sheaf on $X$.
Define the \textit{Kodaira-Iitaka dimension}, $\ka(X,\sL)$, as follows,
\begin{enumerate}
\item
If $h^0(X,\sL^{[\, t]})=0$ for all $t>0$, let $\ka(X,\sL):=-\infty$.
\item
Otherwise, let,
\[
\ka(X,\sL) := \max_{t>0} \dim \phz_{\sL^{[\, t]}}(X), \]

\end{enumerate}
\end{defn}

The Kodaira-Iitaka dimension can also be expressed in terms of the 
assymptotic growth rate of the dimension of spaces of global sections;
at least when we are working in characteristic $0$.
This is well known in the setting of an invertible sheaf on a complete normal variety,
but requires some extra argument for the general situation in definition \ref{KID def}.
\begin{prop}\label{KID assym}
Let $X$ be an almost complete normal variety in charcteristic $0$ 
and $\sL$ a rank one reflexive sheaf on $X$.
If $\ka=\ka(X,\sL)\geq 0$, then there exists constants $0<\al<\be$ for which,
\[
\al t^\ka < h^0(X,\sL^{[\,t]})< \be t^\ka \]
for $t$ sufficiently large and divisible.
\end{prop}
\begin{proof}
We will begin with the case when $X$ is complete.
Assume $\ka(X,\sL)\geq 0$.
Then we may choose a global section $s\in H^0(X,\sL^{[\,m]})$ for some $m>0$.
Then $s$ gives a map,
\[
(\sL^*)^{[\,m]}\stackrel{\ten s}{\>}(\sL^*)^{[\,m]}\ten\sL^{[\,m]}
\> (\sL^*)^{[\,m]}[\ten]\sL^{[\,m]}\iso\sO_X \]
Since $(\sL^*)^{[\,m]}\ten\sL^{[\,m]}$ is torsion-free, this map is injective.
Let $\sI\sub\sO_X$ be the image of $(\sL^*)^{[\,m]}$ in $\sO_X$.
It is not hard to see that for all $t>0$, $s^{[\,t]}\in H^0(X,\sL^{[\,tm]})$
gives a map, $(\sL^*)^{[\,tm]}\stackrel{\sim}{\>}\sI^t\sub\sO_X$.\\
\indent
Let $\pi:\tilde{X}\>X$ be a log resolution with respect to the ideal sheaf $\sI$,
so that $\sI\cdot\sO_{\tilde{X}}\iso\sO_{\tilde{X}}(-D)$ 
for some Cartier divisor $D$ on $\tilde{X}$.
For $t>0$, we will also have that $\sI^t\cdot\sO_{\tilde{X}}\iso\sO(-tD)$.
The kernal of the natural surjection $\pi^*\sI^t\epito\sI^t\cdot\sO_{\tilde{X}}$
will have support contained in the exceptional locus of $\pi$
and thus will be torsion.
This means,
\[
\Hom(\sI^t\cdot\sO_{\tilde{X}},\sO_{\tilde{X}})\iso \Hom(\pi^*\sI^t, \sO_{\tilde{X}}) \]
This allows us to write the following series of isomorphisms,
\begin{align*}
H^0(\tilde{X},\sO_{\tilde{X}}(tD))
&\iso \Hom(\sI^t\cdot\sO_{\tilde{X}},\sO_{\tilde{X}})\\
&\iso \Hom(\pi^*\sI^t, \sO_{\tilde{X}}) \\
&\iso \Hom(\sI^t,\sO_X)\\
&\iso \Hom((\sL^*)^{[\,tm]},\sO_X) \\
&\iso H^0(X,\sL^{[\,tm]}) \end{align*}
Once we have this the result follows from the analogous result
comparing the assymptotic growth rate of $H^0(X,\sO(tD))$
to $\ka(X,\sO(D))$ for a Cartier divisor on a complete smooth variety.
See e.g. \cite{Iitaka}[\S 10.2].\\
\indent
In the general case, let the inclusion $i:X\monoto\overline{X}$
make $X$ into an almost complete normal variety.
In particular, we take $\overline{X}$ to be a complete normal variety.
Then $i_*\sL$ is a rank one reflexive sheaf on $\overline{X}$
and we have equalities,
\begin{align*}
\ka(X,\sL)&=\ka(\overline{X},i_*\sL) \\ 
h^0(X,\sL^{[\,t]})&=h^0(\overline{X}, (i_*\sL)^{[\,t]}) \quad {\rm for\ all}\ t>0 \end{align*}
Then the result in the general case follows from the case when $X$ is complete.
\end{proof}

For a normal variety $X$ with smooth locus $X_{\sm}\sub X$,
we may observe that the Weil divisors on $X$ restricted to $X_{\sm}$ 
are exactly the Cartier divisors on $X_{\sm}$,
and the rank one reflexive sheaves on $X$ restricted to $X_{\sm}$
are exactly the invertible sheaves on $X_{\sm}$.
Furthermore in each case the objects on $X$ are uniquely determined 
by their restriction to $X_{\sm}$.
Thus by the natural correspondence between Cartier divisors and invertible sheaves,
we see that there is a natural isomorphism $\Cl(X)\iso \RPic(X)$ between the divisor class group of $X$
and the group of rank one reflexive sheaves on $X$  up to isomorphism with reflexive tensor product.
Furthermore there is a natural correspondence between effective Weil divisors on $X$
and global sections, $s\in H^0(X,\sL)$,
of the corresponding rank one reflexive sheaf.
In this paper we will write $\sO_X(D)$ 
for the reflexive sheaf corresponding to a Weil divisor D.\\
\indent
Given a rank one reflexive sheaf $\sL$ on a normal variety $X$
and a linear subspace $L\sub H^0(X,\sL)$
we will define the base locus of $L$, $\bl(L)$,
to be the set-theoretic intersection of all the effective Weil divisors corresponding to 
non-zero sections in $L$ if $L\neq \{0\}$.
If $L=\{0\}$ we will let $\bl(L)$ be all of $X$.
Notice that if $L$ is finite dimensional 
then $\phz_L$ will be defined on the open set $X_{\sm}\setminus\bl(L)$.
We will define the base locus of $\sL$ to be,
\[
\bl(\sL):=\bl(H^0(X,\sL)) \]
We will define the stable base locus of $\sL$ to be,
\[
\sbl(\sL):=\bigcap_{t>0} \bl(\sL^{[\,t]}) \]
Similarly we will say a Weil divisor $F$ is a fixed component of 
a linear subspace $L\sub H^0(X,\sL)$ if $F\leq D$
for every effective Weil divisor $D$ corresponding to a non-zero sections in $L$.
This is the same as
$F\cap X_{\sm}$ being a fixed component of $L\vert_{X_{\sm}}$.
\begin{prop}\label{fix comp}
Let $X$ be a normal variety with rank one reflexive sheaf $\sL$
and $L\sub H^0(X,\sL)$ a finite dimensional linear subspace.
If $F$ is a maximal fixed component of $L$, then,
\[
\phz_L^*\sO(1)=\sL(-F)\]
\end{prop}
\begin{proof}
Since we may work over the smooth locus of $X$,
where every divisor is Cartier,
we may assume that $X$ is smooth.
Let $s_F\in H^0(X,\sO_X(F))$
be the section corresponding to the divisor $F$.
Tensor product by $s_F$ defines a map,
\[
\ten s_F: H^0(X,\sL(-F))\>H^0(X,\sL) \]
The map must be injective and since $F$ 
is a fixed component of $L$, this subspace will be in the image of the map. 
Then $(\ten s_F)^{-1}(L)\sub H^0(X,\sL(-F))$
will be a linear subspace without fixed component which gives the same rational map as $L$.
This means $(\ten s_F)^{-1}(L)$ is a subspace of $\phz_L^*\sO(1)$,
which implies our statement. 
\end{proof}

\subsection{Generalized Kodaira-Iitaka Dimension}\ \\
\indent
The definition of Kodaira-Iitaka dimension in subsection \ref{RS}
is essentially geometric in that it is based on the dimensions of the images of certain rational morphisms.
However proposition \ref{KID assym} shows that Kodiara-Iitaka dimension
is also related to the rate of growth of the space of global sections of certain reflexive sheaves.
These two ways of measuring Kodaira-Iitaka dimension leads one to
two different possible ways of generalizing the invariant to schemes which are not almost complete normal varieties.
In this subsection these generalizations will be worked out in greater detail,
being referred to as geometric and algebraic Kodaira-Iitaka dimension.
These notions will be used in various places throughout the remainder of the paper.

\begin{defn}[Geometric Kodaira-Iitaka Dimension] \label{gKID def}
Let $X$ be a variety
and $\sL$ a rank one reflexive sheaf on $X$.
Assume either that $X$ is normal or that $\sL$ is an invertible sheaf.
Define the \textit{geometric Kodaira-Iitaka dimension}, $\ka_g(X,\sL)$, as follows,
\begin{enumerate}
\item
If $h^0(X,\sL^{[\, t]})=0$
for all $t>0$, let $\ka_g(X,\sL):=-\infty$.
\item
Otherwise, let,
\[
\ka_g(X,\sL):=\max_L\dim \phz_L(X) \]
where $L$ is allowed to range over all non-trivial
finite dimensional linear subspaces of $H^0(X,\sL^t)$ for some $t>0$,
and $\phz_L:X\ratmap\P^N$ is the corresponding rational map.
\end{enumerate}
\end{defn}

\begin{note}
This definition agrees with definition \ref{KID def}
when $X$ is an almost complete normal variety.
\end{note}

Let $X$ and $\sL$ be as in definition \ref{gKID def}
and let $L\sub H^0(X,\sL)$ be a finite dimensional linear subspace
with corresponding rational morphism $\phz_L:X\ratmap\P^N$.
Let $A\sub X$ be a subvariety.
If $\sL$ is not invertible, assume that $A$ is normal and not contained in $\Sing X$.
Also assume that $A$ is not contained in $\bl(L)$.
Then $\phz_L$ will be defined on a dense open subset of $A$.
This means we may speak of the image,
$\phz_L(A)\sub\P^N$.
Furthermore, the rational map $\phz_L\vert_A$
will be the same as the one corresponding to 
the subspace $r(L)\sub H^0(A,\sL\vert_A^{**})$,
where $r:H^0(X,\sL)\>H^0(A,\sL\vert_A^{**})$ is the natural reflexive restriction map.
These observations lead to the following proposition.

\begin{prop}\label{kod im 2}
Let $X,\sL,L,\phz_L, A$ be as above. Then,
\[
\dim \phz_L(A)\leq \ka_g(A,\sL\vert_A^{**}) \]
\end{prop}

We will now consider the second generalization of Kodaira-Iitaka dimension,
some of its properties, and its relation to geometric Kodaira-Iitaka dimension.

\begin{defn}[Algebraic Kodaira-Iitaka Dimension]
Let $X$ be a separated scheme of finite type over a field $k$,
not necessarily algebraically closed, and $\sF$ a coherent reflexive sheaf on $X$.
Assume either that $X$ is integral or that $\sF$ is an invertible sheaf.
In this situation define the \textit{algebraic Kodaira-Iitaka dimension}, $\ka_a(X,\sF)$, as follows,
\begin{enumerate}
\item
If
\[
h^0(X,\sym^{[\, t]}\sF)=\dim_k H^0(X,\sym^{[\, t]}\sF)=0 \]
for all $t>0$, let $\ka_a(X,\sF):=-\infty$.
\item
Otherwise, let $\ka_a(X,\sF)$ be the largest integer $\ka$
such that there exists a constant $0<\al$ with,
\[
\al t^{\ka} < h^0(X,\sym^{[\, t]}\sF) \]
for all sufficiently large and divisible $t$.
\item
If such an $\al$ exists for every integer $\ka$, let $\ka_a(X,\sF):=\infty$.
\end{enumerate}
\end{defn}

\begin{prop}\label{aKID usc}
Let $f:X\> S$ be a flat projective morphism of noetherian schemes,
and let $\sL$ be an invertible sheaf on $X$,
then $\ka_a(X_s,\sL_s)$
is upper semicontinuous as a function on $S$.
\end{prop}
\begin{proof}
For any $t>0$ \cite{Hart1}[III.12.8] tells us that,
\[
h^0(s,\sL^t)=\dim_{k(s)}H^0(X_s,\sL_s^t) \]
is upper semicontinuous.
Since $\ka_a(X_s,\sL_s)$ measures the assymptotic growth of this invariant 
for large and divisible $t$,
it is easy to see that it is upper semicontinuous also.
\end{proof}

\begin{prop}\label{fe aKID}
Let $X,k,\sF$ as in the definition of algebraic Kodaira-Iitaka dimension.
Let $k\sub K$ be a field extension 
and $X_K,\ \sF_K$ the appropriate scheme and sheaf after base extension.
Then,
\[
\ka_a(X_K,\sF_K)=\ka_a(X,\sF) \]
\end{prop}
\begin{proof}
Let $f:X_K\>X$ be the (flat) extension map.
Then if $X$ is integral, $X_K$ will be also.
In this case, since for any coherent sheaf, $\sF'$ on $X$,
\[
\sym^t(f^*\sF')\iso f^*(\sym^t\sF') \quad {\rm and} 
\quad (f^*\sF')^*\iso f^*(\sF'^*) \]
we have,
\[
\sym^{[\,t]}\sF_K\iso f^*(\sym^{[\,t]}\sF) \]
This also holds in general when $\sF$ is an invertible sheaf.\\
\indent
Since $K$ is flat over $k$, we know from \cite[III.9.3]{Hart1} that,
\[
H^0(X_K,\sym^{[\,t]}\sF_K)\iso H^0(X_K,f^*(\sym^{[\,t]}\sF))\iso H^0(X,\sym^{[\,t]}\sF)\otimes_k K \]
and the proposition follows.
\end{proof}

\begin{prop}\label{gKID leq aKID}
Let $X$ and $\sL$ be as in the definition of geometric Kodaira-Iitaka dimension. Then,
\[
\ka_g(X,\sL)\leq \ka_a(X,\sL)\]
\end{prop}
\begin{proof}
If $X$ is normal and $\sL$ is a rank one reflexive sheaf,
then we have,
\begin{align*}
\ka_a(X_{\sm},\sL_{\sm})&=\ka_a(X,\sL) \\
\ka_g(X_{\sm},\sL_{\sm})&=\ka_g(X,\sL) \end{align*}
Thus we may assume that $\sL$ is an invertible sheaf.\\
\indent
Assume $\ka_g(X,\sL)\geq 0$.
Choose $m>0$, $L\sub H^0(X,\sL^m)$
so that $\ka_g(X,\sL)=\dim \phz_L(X)$.
By the definition of $\phz_L$ we have,
\[
H^0(\phz_L(X),\sO(t))\iso L^t\sub H^0(X,\sL^{tm}) \]
Since the Hilbert polynomial of $\phz_L(X)$ has degree $\dim \phz_L(X)=\ka_g(X,\sL)=\ka_g$, 
we have for some $\al>0$ and all $t\gg 0$,
\begin{align*}
\al t^{\ka_g} & < h^0(\phz_L(X),\sO(t)) \\
&\leq h^0(X,\sL^{tm}) \end{align*}
From this the proposition follows.
\end{proof}

\begin{prop}\label{KID=gKID=aKID}
If $X$ is an almost complete normal variety in characteristic $0$
and $\sL$ is a rank one reflexive sheaf on $X$, then,
\[
\ka(X,\sL)=\ka_g(X,\sL)=\ka_a(X,\sL) \]
\end{prop}
\begin{proof}
This follows from the discussion so far and in particular proposition \ref{KID assym}.
\end{proof}

\section{Kodaira-Iitaka Dimension on a Subvariety}\label{Part I}
\subsection{Divisorial Valuations}\ \\
\indent
Next we will consider prime divisors on a normal variety and 
also divisorial valuations.
A prime divisor, that is, a closed integral subvariety of codimension one, on a normal variety
determines a discrete valuation ring with the function field of the variety
as its quotient field.
This may be abstracted to define the concept of divisorial valuation rings.

\begin{defn}
Let $K$ be a finitely generated field 
over an algebraically closed field $k$ with $\ptd_k K =n$,
let $\nu$ be a discrete valuation on $K$,
let $R\sub K$ be the corresponding valuation ring,
and let $m_R$ be the maximal ideal of $R$. If,
\[
\ptd_k R/m_R = n - 1 \]
we will say that $\nu$ is a \textit{divisorial valuation} on $K$
and that $R$ is a \textit{divisorial valuation ring} with quotient field $K$.\\
\indent
If $X$ is a variety with function field $K$ and
$Y$ a prime divisor on $X$ such that $(\sO_{X, Y},m_Y)\iso(R,m_R)$ in $K=K(X)$
we will say that $(X,Y)$ is a \textit{model} for $\nu$ or for $R$.
\end{defn}

It is not a priori clear that a given divisorial valuation has a model.
However if $X$ is a complete variety and $R$ is a divisorial valuation ring with quotient field $K(X)$
we may always associate to $R$ a closed subvariety of $X$.
In such a case $X$ is proper over $\Spec k$.
Since there is natural map $\Spec K(X)\> X$,
we know by the valuative criterion of properness that $\Spec R\>\Spec k$
has a unique lift, $f:\Spec R\> X$.

\begin{defn}
We will give the name $\cen_X R$ to the closure of the image of the unique closed point in $\Spec R$
under the above morphism, $f:\Spec R\> X$.
Then $\cen_X R$ will be a closed subvariety in $X$.
\end{defn}

The following proposition is Lemma 2.45 in \cite{K+M}.
\begin{prop}\label{res ad}
Let $X$ be a complete variety over $k$,
and $R$ be a divisorial valuation ring with quotient field $K(X)$.
Define inductively,
\begin{enumerate}
\item
$X_0:=X$
\item
If $(X_i, \cen_{X_i} R)$ is not a model for $R$, let $\pi_{i+1}:X_{i+1}\>X_i$
be the blow-up of $X_i$ along $\cen_{X_i} R$
\end{enumerate}
This process eventually terminates.
That is, for some $r\geq 0$ $(X_r, \cen_{X_r} R)$ is a model for $R$.
\end{prop}

We will now consider the restriction of a divisorial valuation ring to a subfield
of its quotient field.
This will be useful for studying the image
of a prime divisor under a dominant rational map.\\
\indent
Let $k$ be an algebraically closed field and let $(K,R)$ be a pair, where
$K$ is a finitely generated field over $k$ with $\ptd_k K = n$
and $R$ is a divisorial valuation ring with quotient field $K$
and maximal ideal $m_R$, corresponding to a valuation $\nu$ on $K$.
Consider a subfield, $K'\sub K$
with $\ptd_k K' = m\leq n$.
It is possible that $m_R\cap K'=\{0\}$.
If this is not the case, $\nu\vert_{K'}$ is a valuation on $K'$.
Then $R\cap K'\sub K'$ is a discrete valuation ring
with maximal ideal $m_R\cap K'$.
\begin{prop}
Let $K,R,\nu,K',m_R$ be as above.
If $m_R\cap K'\neq\{0\}$, then $\nu\vert_{K'}$ is a divisorial valuation on $K'$
with valuation ring $R\cap K'$.
\end{prop} 
\begin{proof}
We only need to show that $\ptd_k (R\cap K')/(m_R\cap K') = m-1$.
We easily see that $\ptd_k (R\cap K')/(m_R\cap K') \leq m-1$,
since $m_R\cap K' \neq \{0\}$ is a prime ideal of height at least $1$ in $R\cap K'$.\\
\indent
Let $\{f_1,\ldots, f_r\}\sub R\cap K'$ be chosen
such that $\{\overline{f_i}\}$ is a trancendence basis for 
$(R\cap K')/(m_R\cap K')\sub R/m_R$ over $k$.
Then complete a trancendence basis for $R/m_R$ by choosing
$\{g_1,\ldots, g_s\}\sub R$ so that $\{\overline{f_i}, \overline{g_j}\}$
is a trancendence basis for $R/m_R$ over $k$.
\begin{lem}
In the above, the elements of $\{g_j\}$ are algebraically independent and trancendent over $K'$.
\end{lem}
\begin{proof}
Suppose otherwise, then there would be a nontrivial polynomial
\[
F(x_1,\ldots x_s)=\sum c_ix_1^{a_{i,1}}\cdots x_s^{a_{i,s}} \]
with coefficients in $K'$ and $F(g_1,\ldots, g_s)=0$.
Consider the valuations, $\nu(c_i)$, of the coefficients of $F$.
Let the minimum be achieved by $c_{\min}$.
Then $c_{\min}^{-1}F$ has coefficients in $R\cap K'$,
and $\overline{c_{\min}^{-1}F}$ is a nontrivial polynomial
with coefficients in $(R\cap K')/(m_R\cap K')$ such that
$F(\overline{g_1},\ldots,\overline{g_s})=0$.
This would make the elements of $\{\overline{g_j}\}$ 
algebraically dependent over $(R\cap K')/(m_R\cap K')$
and thus over $k(\overline{f_i})$.
But this contradicts the choice of $\{g_j\}$.
\end{proof}
The lemma implies 
\[
s\leq \ptd_{K'} K = n-m \]
Also we know 
\[
r+s=n-1=\ptd_k R/m_R \]
Thus we have $r\geq m-1$,
proving the proposition.
\end{proof}
We wish to consider the geometric consequences of this algebraic fact.
We immediately have the following.

\begin{cor}\label{rat ad}
Let $X$ and $V$ be varieties and
$\phz:X\ratmap V$ a dominant rational map
corresponding to $K(V)\sub K(X)$,
let $R$ be a divisorial valuation ring with quotient field $K(X)$.
Then either $R\cap K(V)=K(V)$ or $R\cap K(V)$ is
a divisorial valuation ring with quotient field $K(V)$.
\end{cor}

If in the above situation $R$ has a model of the form $(X,Y)$,
we may reach the following conclusion

\begin{prop}\label{rat im Y}
Let $X, V, \phz, R$ be as in Corollary \ref{rat ad}.
Suppose that $(X,Y)$ is a model for $R$ and that $V$ is complete.
If $R\cap K(V)=K(V)$ then $\phz(Y)=V$.
If $R\cap K(V)\neq K(V)$ then
$\phz(Y)=\cen_V R\cap K(V)$.
\end{prop}

\begin{proof}
If $R\cap K(V)=K(V)$ then the generic point of $Y$ is taken to the generic point of $V$;
in other words, $\phz(Y)=V$.
If $R\cap K(V)\neq K(V)$ then $R\cap K(V)$ is a divisorial valuation ring with
quotient field $K(V)$.
It is not hard to see that $\phz(Y)=\cen_V R\cap K(V)$.
\end{proof}

\subsection{Kodaira-Iitaka Dimension on a Prime Divisor}
\begin{prop}\label{main kod}
Let $X$ be a normal variety and
$Y\sub X$ a prime divisor.
Let $\phz:X\ratmap\P^N$ be a rational map.
Then, either $\phz(Y)=\phz(X)$,
or we may choose $\epz>0$ so that,
\[
\dim \phz(X) -1 \leq\ka_g(Y, (\phz^*\sO_{\P^N}(n)[\otimes]\sO_X(-bY))\vert_Y^{**})
\quad \text{for\ all}\ n,b > 0,\ \frac{b}{n}<\epz\]
\end{prop}
\begin{proof}
With the situation as above.
Let $V = \phz(X)\sub\P^N$.
Then $V$ is complete
and the rational map $\phz$ corresponds to the inclusion $K(V)\sub K(X)$.
Let $R=\sO_{X,Y}\sub K(X)$.
To prove the theorem we may assume $\dim \phz(Y)<\dim V$,
since otherwise we are done.\\
\indent
First suppose that $\phz(Y)\sub V$ is a divisor.
Let $H$ be the divisor on $V$ corresponding to $\sO_V(1)$,
and $Z$ be a Cartier divisor on $V$ with $Z\geq \phz(Y)$.
Since $V$ is projective, its cone of ample divisors is open in $\N1(V)$.
This means we may choose $\epz>0$ so that 
$H-\frac{b}{n}Z$ is ample as a $\Q$-divisor for all $n,b>0,\ \frac{b}{n}<\epz$.
Thus the line bundle $(\sO_V(n)\otimes\sO_V(-bZ))\vert_{\phz(Y)}$ is ample.
It follows that,
\begin{align*}
\dim \phz(X)-1 &=\dim \phz(Y)\\
&\leq \ka_g(Y,\phz^*((\sO_V(n)\otimes\sO_V(-bZ))\vert_{\phz(Y)}))\\
&\leq \ka_g(Y,(\phz^*\sO_{\P^N}(n)[\otimes]\sO_X(-bY))\vert_Y^{**})
\end{align*}

\indent
Assume then that $\phz(Y)$ is neither $V$ nor a divisor on $V$.
By Corollary \ref{rat ad} and Proposition \ref{rat im Y}, 
$R\cap K(V)$ is a divisorial valuation ring with quotient field $K(V)$,
which has center $\phz(Y)$.
Then according to Proposition \ref{res ad}, there is a resolution by blow-ups,
\[
\tilde{V}=V_r\stackrel{\pi_r}{\>}V_{r-1}\>\cdots\>V_1\stackrel{\pi_1}{\>}V_0=V \]
such that $(\tilde{V}, \cen_{\tilde{V}}R\cap K(V))$ is a model for $R\cap K(V)$.
Let $W_i\sub V_i$ be the exceptional divisor of the blow-up $\pi_i$.
By a slight abuse of notation, let $W_i$ also represent the pull-back of $W_i$ to $\tilde{V}$.
Note that these $W_i$ are Cartier.

\begin{lem}\label{amp tildeV}
Let $H$ be a very ample divisor on $V$,
then for some $\{\al_i\in\Q\},\ 1\gg\al_1\gg\ldots\gg\al_r>0$ we have that,
\[
\pi^*H-\sum\al_i W_i \]
is ample as a $\Q$-divisor on $\tilde{V}=V_r$.
\end{lem}

\begin{proof}
According to \cite[II.7.10]{Hart1} and the fact that $\sO_{V_1}(1)=\sO_{V_1}(-W_1)$,
we have that $n_1\pi_1^*H-W_1$ is very ample on $V_1$ for $n_1\gg 0$.
Thus $\pi_1^*H-\al_1 W_1$ is ample, where we can choose $\al_1\in\Q$ arbitrarily small and positive.
Continuing in this manner we get our result, where
at each step we can choose $\al_i$ arbitrarily small compared to the $\al$'s already chosen.
\end{proof}

We may however say more than this.
In particular, let us consider the set of possible $\{\al_i\}$
for which $\pi^*H-\sum\al_i W_i$ is ample as a $\Q$-divisor.
There is a natural linear map $f:\R^r\>\N1(\tilde{V})$,
\[
f:(\al_1,\ldots, \al_r) \mapsto \pi^*H-\sum\al_i W_i \]
Since $\tilde{V}$ is projective,
its cone of ample divisors, $\Amp(\tilde{V})\sub\N1(\tilde{V})$,
is open and convex.
This means the inverse image $f^{-1}(\Amp(\tilde{V}))\sub\R^r$
is also open and convex.
Also, according to lemma \ref{amp tildeV},
we can choose $1\gg\al_1\gg\ldots\gg\al_r>0$
so that $(\al_1,\ldots,\al_r)\in f^{-1}(\Amp(\tilde{V}))$.\\
\indent
Now let $\tilde{\phz}: X\ratmap\tilde{V}$ be the appropriate rational map.
We know from proposition \ref{rat im Y} that
$\tilde{\phz}(Y)=\cen_{\tilde{V}}R\cap K(V)$.
Since this is a model for an abstract prime divisor,
we know $\dim \tilde{\phz}(Y)=\dim \tilde{V}-1=\dim V -1$.\\
\indent
Since $V\sub\P^N$, $\sO_{\P^N}(1)\vert_V=\sO_V(1)$ is a very ample divisor on $V$.
Let $L$ be the divisor on $X$ corresponding to $\phz^*\sO_{\P^N}(1)$.
As in the discussion above, choose $1\gg\al_1\gg\ldots\gg\al_r>0$ and $n \gg 0$ so that
$(\pi^*\sO_V(n)\otimes\sO_{\tilde{V}}(-n\sum\al_i W_i))$ is very ample on $\tilde{V}$.
Let $\psi:\tilde{V}\>\P^M$ be the embedding corresponding to this divisor.
\[
\xymatrix{
&\tilde{V}\ar[r]^\psi\ar[d]^\pi&\P^M\\
X \ar@{-->}[ur]^{\tilde{\phz}}\ar@{-->}[r]_\phz & V & }
\]
Then $\psi\circ\tilde{\phz}:X\ratmap \P^M$ is a rational map
with $\dim \psi\circ\tilde{\phz}(Y)=\dim V-1$.
Note that $\tilde{\phz}^*\pi^*\sO_V(1)=\phz^*\sO_V(1)$,
since they agree on the domain of definition of $\tilde{\phz}$, $U$.
Therefore we have,
\[
(\psi\circ\tilde{\phz})^*\sO_{\P^M}(1) = \sO_X(n(L-\sum\al_i\tilde{\phz}^*W_i))\]
For each $W_i$ let $\tilde{\phz}^*W_i = D_i+m_iY$,
where $D_i$ is an effective (possibly zero) divisor on $X$ not containing $Y$ in its support, $m_i\geq 0$ and $m_r>0$
since $\tilde{\phz}(Y)\sub W_r$.
Then $\psi\circ\tilde{\phz}$ is the rational map corresponding to a linear subspace of,
\[
H^0\left(X, \sO_X\left(n\left(L-\sum\al_iD_i-\left(\sum \al_i m_i\right) Y\right)\right)\right) \]
Since the rank one reflexive sheaf above is the pull back of $\sO_{\P^M}(1)$, 
it will not have any fixed components
and $Y$, being a divisor, will not be contained in its base locus.
Neither will $Y$ be contained in $\Sing X$ which has codimension at least $2$.
Therefore by proposition \ref{kod im 2}, we have,
\[
\dim V -1 =
\dim \psi\circ\tilde{\phz}(Y) \leq
\ka_g(Y,\sO_X\left(n\left(L-\sum\al_iD_i-\left(\sum \al_i m_i\right) Y\right)\right)\big\vert_Y^{**})\]
Also each $D_i$ intersects $Y$ non-negatively,
so we have,
\begin{align*}
&\ka_g(Y,\sO_X\left(n\left(L-\sum\al_iD_i-\left(\sum \al_i m_i\right) Y\right)\right)\big\vert_Y^{**})
\\
\leq\ &\ka_g(Y, \sO_X\left(n\left(L-\left(\sum\al_i m_i \right) Y\right)\right)\big\vert_Y^{**})
\end{align*}
\indent
We are interested in the possible values of $\sum\al_i m_i$.
This quantity may be chosen to be any value in the image of $f^{-1}(\Amp(\tilde{V}))$, as above,
under the linear map,
\[
(\al_1,\ldots,\al_r)\mapsto \sum \al_i m_i \]
Since such a linear projection is an open map,
the image will be open and convex in $\R$.
Also, according to lemma \ref{amp tildeV}, it contains arbitrarily small positive values.
Thus there exists $\epz>0$ so that $\sum \al_i m_i$ may be taken to be any
rational number in the open interval $(0,\epz)$.
In terms of rank one reflexive sheaves this means,
\[
 \dim \phz(X) -1 \leq
\ka_g(Y, \left(\phz^*\sO_{\P^N}(n)[\otimes]\sO_X(-bY)\right)\vert_Y^{**})
\quad \text{for all}\ n,b>0,\ \frac{b}{n}<\epz \]
\end{proof}

\begin{thm}\label{main div}
Let $X$ and $Y\sub X$ be as in Proposition \ref{main kod} and
let $\sL$ be a rank one reflexive sheaf on $X$.
If $\ka_g(X,\sL)\geq 0$, then for $a>0$, sufficiently large and divisible, there exists an integer $m\geq 0$,
which is taken to be zero if $Y$ is not a fixed component of $\sL^{[\,a]}$,
and real number $\epz>0$
such that one of the following holds,
\[
\ka_g(X,\sL) \leq
\ka_g(Y, (\sL^{[\,a]}[\otimes]\sO(-Y)^{[\,m]})\vert_Y^{**}) \]
or,
\[
0 \leq
\ka_g(Y, (\sL^{[\,a]}[\otimes]\sO(-Y)^{[\,m]})\vert_Y^{**})\ \text{and,}\]
\[
\ka_g(X,\sL)-1 \leq
\ka_g(Y,(\sL^{[\,na]}[\otimes]\sO_X(-Y)^{[\,nm+b]})\vert_Y^{**})
\quad \text{for\ all} \ n,b>0,\ \frac{b}{n}<\epz \]
\end{thm}

\begin{proof}
If $\ka_g(X,\sL)=-\infty$ we are done.
Thus we may assume that for $a$ large and divisible,
there is a subspace $L\sub H^0(X,\sL^{[\, a]})$
such that $\phz_L:X\ratmap\P^N$ is a rational map
with $\ka_g(X,\sL)=\dim \phz_L(X)$.
Let $L$ have maximum fixed component $F=D+mY$.
If $Y$ is not a fixed component of $\sL^{[\, a]}$,
we can choose $L$ so that $m=0$.
By proposition \ref{fix comp}
we know that $\sL^{[\,a]}[\otimes]\sO_X(-F)=\phz_L^*\sO_{\P^N}(1)$.\\
\indent
We have by proposition \ref{kod im 2} and the fact that $h^0(Y,\sO_X(D)\vert_Y^{**})>0$,
\begin{align*}
\dim \phz_L(Y)
&\leq \ka_g(Y, (\sL^{[\,a]}[\otimes]\sO_X(-F))\vert_Y^{**})\\
&\leq \ka_g(Y,(\sL^{[\,a]}[\otimes]\sO_X(-Y)^{[\,m]})\vert_Y^{**})
\end{align*}
If $\phz_L(Y)=\phz_L(X)$ we have,
\[
\ka_g(X,\sL)\leq\ka_g(Y, (\sL^{[\,a]}[\otimes]\sO(-Y)^{[\,m]})\vert_Y^{**}) \]
Otherwise, we still have,
\[
0\leq\ka_g(Y, (\sL^{[\,a]}[\otimes]\sO(-Y)^{[\,m]})\vert_Y^{**})\]
Also, by Proposition \ref{main kod}, in this case,
\[
\ka_g(X,\sL)-1=\dim\phz(X)-1\leq
\ka_g(Y, \left((\sL^{[\,a]}[\otimes]\sO_X(-F))^{[\,n]}[\otimes]\sO(-bY)\right)\vert_Y^{**})\]
for all $n,b>0,\ \frac{b}{n}<\epz$.
However,
\begin{align*}
&\ka_g(Y, \left((\sL^{[\,a]}[\otimes]\sO_X(-F))^{[\,n]}[\otimes]\sO(-bY)\right)\vert_Y^{**})\\
=\ &\ka_g(Y,(\sL^{[\,na]}\otimes\sO_X(-nD-(nm+b)Y))\vert_Y^{**}) \\
\leq\ &\ka_g(Y,(\sL^{na}\otimes\sO_X(-Y)^{[\,nm+b]})\vert_Y^{**}) 
\end{align*}

\end{proof} 

\subsection{Application to Subvarieties of any Codimension}\label{apps}\ \\
\indent
By using blow-ups, we may apply these results to the case of a small normal pair $(X,A)$.
We will start with the case when both $X$ and $A$ are smooth, $A$ is closed in $X$, and $\sL$ is an invertible sheaf.
In order to avoid an overly cumbersome statement, the result will be simplified and slightly weakened here,
but it will be easy to see how Theorem \ref{main div}
yields a slightly stronger statement than given here.

\begin{prop}\label{smooth subvar}
If $X$ is a smooth variety,
$A\sub X$ a smooth closed subvariety of codimension $d$
with normal sheaf $\sN_{A\vert X}$, and
$\sL$ an invertible sheaf on $X$,
then for some integers $n_1>0,\ n_2\geq 0$,
\[
\ka_g(X,\sL)-1\leq
\ka_a(A,\sL^{n_1}\vert_A\ten\sym^{n_2}\sN_{A\vert X}^*) \]
Furthermore, if $A$ is not contained in the stable base locus of $\sL$, we may choose $n_1\gg n_2$.
\end{prop}
\begin{proof}
In the case $d=1$ this follows immediately from theorem $\ref{main div}$
and proposition \ref{gKID leq aKID},
taking $n_1=na,\ n_2= nm+b$ or $nm$;
thus assume $d>1$.\\
\indent
Let $\pi:\tilde{X}\>X$ be the blow-up of $X$ along $A$.
Let $Y\sub\tilde{X}$ be the exceptional divisor.
Then according to \cite[II.8.24]{Hart1}, $\pi:Y\>A$
is isomorphic to the projective space bundle, $\P(\sN_{A\vert X}^*)$
with $\sN_{Y\vert\tilde{X}}^*$ corresponding to $\sO_{\P(\sN_{A\vert X}^*)}(1)$.
In particular this means by \cite[II.7.11]{Hart1} that,
\[
\pi_*\sN_{Y\vert\tilde{X}}^{-k}\iso\sym^k\sN_{A\vert X}^* \]
Since $\pi_*\sO_{\tilde{X}}=\sO_X$, we have,
$\ka_g(X,\sL)=\ka_g(\tilde{X},\pi^*\sL)=\ka$.\\
\indent 
We may apply theorem \ref{main div} to conclude,
\[
\ka-1=\ka_g(\tilde{X},\pi^*\sL)-1 \leq
\ka_g(Y,\pi^*\sL^{n_1}\vert_Y\otimes\sN_{Y\vert\tilde{X}}^{-n_2})\]
for $n_1=na$, $n_2=nm+b$ or $nm$. \\
\indent
In terms of the assymptotic growth of global sections, this means,
\[
\al t^{\ka-1}< h^0(Y,\pi^*\sL^{tn_1}\vert_Y\otimes\sN_{Y\vert\tilde{X}}^{-tn_2}) \]
for some $\al>0$ and $t$ large and divisible.\\
\indent
By an application of the projection formula we have,
\begin{align*}
h^0(Y,\pi^*\sL^{tn_1}\vert_Y\otimes\sN_{Y\vert\tilde{X}}^{-tn_2})
&=h^0(A, \pi_*(\pi^*\sL^{tn_1}\vert_Y\otimes\sN_{Y\vert\tilde{X}}^{-tn_2}))\\
&=h^0(A,\sL^{tn_1}\vert_A\otimes\pi_*\sN_{Y\vert\tilde{X}}^{-tn_2})\\
&=h^0(A,\sL^{tn_1}\vert_A\otimes\sym^{tn_2}\sN_{A\vert X}^*)\\
&=h^0(A,\sym^t(\sL^{n_1}\ten\sym^{n_2}\sN_{A\vert X}^*))
\end{align*}
From this the first statement in the proposition follows.\\
\indent
If $A$ is not in the stable base locus of $\sL$, then we may choose $a$
so that $Y$ is not a fixed component of $\pi^*\sL^a$.
Thus in this case we may choose $n_1\gg n_2 = b$ or $0$.
\end{proof}

\begin{thm}\label{main subvar}
Let $(X,A)$ be a small normal pair.
Let $\sN_{A\vert X}$ be the normal sheaf of $A$ in $X$, and
$\sL$ a rank one reflexive sheaf on $X$.
Then for some integers $n_1>0,\ n_2\geq 0$,
\[
\ka(X,\sL)-1 \leq
\ka_a(A,\sL^{[\, n_1]}\vert_A^{**}[\ten]\sym^{[\, n_2]}\sN_{A\vert X}^*)\]
Furthermore, if $A$ is not contained in the stable base locus of $\sL$ then we may choose $n_1\gg n_2$.
\end{thm}

\begin{proof}
Let $U= X\setminus (\Sing X \cup \Sing A\cup\overline{A}\setminus A)$.
Then $U$ and $A\cap U$ are smooth,
$A\cap U$ is closed in $U$,
$\sL\vert_U$ is an invertible sheaf,
$\codim_X(X\setminus U)\geq 2$,
and $\codim_A(A\setminus(A\cap U))\geq 2$.\\
\indent
Then we may apply proposition $\ref{smooth subvar}$ to $A\cap U\sub U$
and $\sL\vert_U$
and observe that,
\[
\ka_a(A, \sL^{[\, n_1]}\vert_A^{**}[\otimes] \sym^{[\, n_2]}\sN_{A\vert X}^*)
=\ka_a(A\cap U, \sL^{n_1}\vert_{A\cap U}\otimes \sym^{n_2}\sN_{A\cap U\vert U}^*) \]
and,
\[
\ka(X,\sL)=\ka(U,\sL\vert_U)=\ka_g(U,\sL\vert_U) \]
to conclude what we desire.

\end{proof}

\subsection{Connections with the Related Results of T. Peternell, M. Schneider and A.J. Sommese}\ \\
\indent
In \cite{PSS}, T. Peternell, M. Schneider and A.J. Sommese
establish an inequality similar to the one in theorem \ref{main subvar} \cite[theorem 2.1]{PSS}.
They use this inequality as a basis for proving several results relating 
Kodaira-Iitaka dimension on a variety to Kodaira-Iitaka dimension on a subvariety.
These results can be proved using theorem \ref{main subvar} in place of \cite[theorem 2.1]{PSS}.\\
\indent
Before stating these results we need to introduce three definitions.
\begin{defn}[$\Q$-Effective]
A coherent sheaf $\sF$ on an almost complete normal variety is said to be $\Q$\textit{-effective}
if these exists an integer $m>0$ such that $\sym^{[\,m]}(\sF^{**})$ is generically spanned by global sections.
\end{defn}
\begin{defn}[Generically Nef]
A coherent sheaf $\sF$ on an almost complete normal variety, $X$ is said to be \textit{generically nef}
if for some open $U\sub X$ with $\codim_X(X\setminus U)\geq 2$, $\sF\vert_U$ is locally free and nef.
\end{defn}
\begin{defn}[Arithmetic Kodaira Dimension]
The \textit{arithmetic Kodaira dimension} of an almost complete normal variety is defined to be,
\[
\ka(X):=\ka(X, \sO(K_X)) \]
\end{defn}
\begin{note}
If $X$ is smooth, $\ka(X)$ agrees with the usual Kodaira dimension.
If $X$ is any normal variety and $\tilde{X}\>X$ is a resolution of singularities,
then,
\[
\ka(\tilde{X})\leq \ka(X) \]
But this inequality may be strict in some circumstances.
\end{note}
\begin{thm}[see \protect{\cite[theorem 4.1]{PSS}}]\label{PSS1}
Let $(X,A)$ be a small normal pair in characteristic $0$
and $\sL$ a rank one reflexive sheaf on $X$. If $\sN_{A\vert X}$ is $\Q$-effective,
then,
\[
\ka(X,\sL)\leq \ka(A,\sL\vert_A^{**})+\codim_X A \]
In particular,
\[
\ka(X)\leq \ka(A)+\codim_X A \]
\end{thm}
\begin{proof}
This is an application of theorem \ref{main subvar}.
Restricting ourselves to characteristic $0$
allows us to apply proposition \ref{KID assym} at the appropriate places.\\
\indent
Let $\ka=\ka(X,\sL)$ and $d=\codim_X A$.
For $t$ large and divisible we have,
\[
\al t^{\ka-1}\leq h^0(A,\sL^{[\,tn_1]}\vert_A^{**}[\otimes]\sym^{[\,tn_2]}\sN_{A\vert X}^*)\]
for some $\al>0$, $n_1>0$ and $n_2\geq 0$.
If $n_2=0$, we have,
\[
\al t^{\ka-1}<h^0(A,\sL^{[\,tn_1]}\vert_A^{**}) \]
This would mean $\ka-1=\ka(X,\sL)-1\leq\ka(A,\sL\vert_A^{**})$,
which establishes the theorem.\\
\indent
Thus assume $n_2>0$.
Since $\sN_{A\vert X}$ is $\Q$-effective, 
we may assume that for $t$ large and divisible enough,
$\sym^{[\,tn_2]}\sN_{A\vert X}$ is generically spanned by global sections.
This means $\sym^{[\,tn_2]}\sN_{A\vert X}^*\monoto\sO_A^{\oplus M}$,
where $M=\rk\sym^{[\,tn_2]}\sN_{A\vert X}<\be_1(tn_2)^{d-1}$.
Then we have,
\begin{align*}
h^0(A,\sL^{[\,tn_1]}\vert_A^{**}[\otimes]\sym^{[\,tn_2]}\sN_{A\vert X}^*)
&\leq h^0(A,\sL^{[\,tn_1]}\vert_A^{**}[\otimes]\sO_A^{\oplus M})\\
&= Mh^0(A,\sL^{[\,tn_1]}\vert_A^{**})\\
&< \be_1(tn_2)^{d-1}\be_2(tn_1)^{\ka(A,\sL\vert_A^{**})}
\end{align*}
Putting this together we have,
\[
\al t^{\ka-1}<\be_1\be_2n_1^{\ka(A,\sL\vert_A^{**})}n_2^{d-1}t^{\ka(A,\sL\vert_A^{**})+d -1} \]
for $t$ large and divisible.
This means, 
\[
\ka \leq \ka(A,\sL\vert_A^{**})+d \]
The second statement in the theorem follows since for $\sN_{A\vert X}$ $\Q$-effective
we have, 
\[
\ka(A,\sO(K_X)\vert_A^{**})\leq \ka(A) \] 
\end{proof}
The next two results follow from theorem \ref{main subvar}, \cite[lemma 3.8]{PSS}, \cite[lemma 3.9]{PSS} and \cite[lemma 3.10]{PSS}.
Their proofs may be found in \cite{PSS}.
In those proofs theorem \ref{main subvar} may be used in place of \cite[theorem 2.1]{PSS}
to produce the same arguments.
\begin{thm}[\protect{\cite[theorem 4.3]{PSS}}]
Let $(X,A)$ be a small normal pair over $\C$.
Suppose $A$ lies as an open subset in a complete normal variety $\overline{A}$
with $\codim_{\overline{A}}(\overline{A}\setminus A)\geq 2$, such that
$\overline{A}$ has at worst terminal singularities and admits a good minimal model.
If $\sN_{A\vert X}$ is generically nef, then,
\[
\ka(X)\leq \ka(A)+\codim_X A \]
\end{thm}
\begin{thm}[\protect{\cite[theorem 4.4]{PSS}}]
Let $(X,A)$ be a small normal pair over $\C$ and $\sL$ a rank one reflexive sheaf on $X$.
Suppose $\sL$ is $\Q$-Cartier and $\sL\vert_A^{**}$ is semi-ample.
If $\sN_{A\vert X}$ is generically nef, then
\[
\ka(\sL)\leq \ka(\sL\vert_A^{**})+\codim_X A \]
\end{thm}

To give more detail we may look at the inequality of interest in \cite{PSS}.
\begin{thm}[\protect{\cite[theorem 2.1]{PSS}}]\label{PSS ineq}
Let $(X, A)$ be a small normal pair over $\C$ and $\sL$ a rank one reflexive sheaf on $X$.
Then there is a positive integer $c$ such that for all $t\geq 0$,
\[
h^0(X,\sL^{[\,t]})\leq \sum_{k=0}^{ct}h^0(A,\sL^{[\,t]}\vert_A^{**}[\otimes]\sym^{[\,k]}\sN_{A\vert X}^*)\]
\end{thm}
A similar inequality can be shown as a corollary of theorem \ref{main div} using the methods in subsection \ref{apps}.
This inequality, which is below, has the advantage that the constant $c$, 
controlling the contribution from $\sN_{A\vert X}^*$ compared to that from $\sL\vert_A^{**}$,
can be chosen to be arbitrarily small
when $A$ is not in the stable base locus of $\sL$.
However it has the disadvantage that we must choose a scaling factor $\be>0$, making the inequality
assymptotic rather than precise.
\begin{cor}[of theorem \ref{main div}]
Let $(X, A)$ be a small normal pair and $\sL$ a rank one reflexive sheaf on $X$.
Then there exists a constant $\be>0$ and positive rational number $c$, which can be chosen to be arbitrarily
small if $A$ is not in the stable base locus of $\sL$,
such that for all sufficiently large and divisible $t$,
\[
h^0(X,\sL^{[\,t]})\leq 
\be \sum_{k=0}^{\flr{ct}}h^0(A,\sL^{[\,t]}\vert_A^{**}[\otimes]\sym^{[\,k]}\sN_{A\vert X}^*)\]
\end{cor}
\begin{proof}
This follows by using the full extent of theorem \ref{main div},
applying the blowing-up methods of subsection \ref{apps}
and doing some arithmetic to fit the results into the form of the sum in the corollary.
\end{proof}

\section{Kodaira-Iitaka Dimension on Intersecting Pairs of Subvarieties}\label{Part II}
\subsection{Inequality for Subvarieties Intersecting Nicely}\ \\
\indent
In this section we will consider the situation
of a pair of subvarieties, $A$ and $B$, in $X$.
We will work over characteristic $0$ throughout this section
in order for morphisms to have dense open smooth loci.
We will begin by finding an inequality relating the dimensions
of the images of $X$, $A$, and $B$ under a rational map
in the case that $A$ and $B$ intersect in a very nice way with respect to this map.\\
\indent
Let $X$ be a variety over an algebraicallly closed field $k$
of characteristic $0$.
Let $\phz:X\ratmap V$ be a dominant rational map, $U\sub X$ the open set on which $\phz$ is defined
and $U_{\sm}\sub U$ a smooth open set
on which $\phz$ is a smooth morphism.
If $A\sub X$ is any subvariety with $A\cap U\neq\emptyset$,
let $\phz(A)\sub V$ be the closure of $\phz\vert_U(A\cap U)$ in $V$.

\begin{defn}[Transverse Intersection]
Let $A,B\sub X$ be subvarieties of a variety $X$.
We will say $A$ and $B$ \textit{intersect transversely} at a closed point, $x\in X$,
if $x\in A\cap B$, and
$T_{X,x}=T_{A,x}+T_{B,x}$.
\end{defn}

\begin{prop}\label{dim V leq}
If $A,B\sub X$ are subvarieties which intersect transversely
at a closed point $x\in U_{\sm}$ such that
$\phz(A)$ and $\phz(B)$ are smooth at $\phz(x)\in V$
then,
\[
\dim V\leq \dim\phz(A)+\dim\phz(B) \]
\end{prop}
\begin{proof}
Since $\phz$ is smooth on $U_{\sm}$,
the tangent map, $T_\phz:T_{X,x}\>T_{V,\phz(x)}$, will be surjective.
By transversality we have,
\[
T_{A,x}+T_{B,x}=T_{X,x} \]
\[
T_\phz T_{A,x}+T_\phz T_{B,x} = T_\phz T_{X,x} = T_{V,\phz(x)} \]
We necessarily have $T_\phz T_{A,x}\sub T_{\phz(A),\phz(x)}$. 
Therefore, since we assume $\phz(A)$ and $\phz(B)$ are smooth at $\phz(x)$, we have,
\begin{align*}
\dim V &\leq \dim T_{V,\phz(x)} \\
&= \dim (T_{\phz(A),\phz(x)}+T_{\phz(B),\phz(x)}) \\
&\leq \dim\phz(A)+\dim\phz(B)
\end{align*}
\end{proof}

\begin{lem}\label{trans cond}
Let $A,B\sub X$ be subvarieties and $Y=A\cap B$ be their scheme theoretic intersection.
Suppose there exists a smooth open subset $W\sub Y$ such that
$\dim W=\dim A+\dim B-\dim X$.
Then $A$ and $B$ intersect transversely at every closed point $x\in W$ which is smooth in $X$.
\end{lem}
\begin{proof}
Let $x\in W$ be a closed point, smooth in $X$.
Suppose $A$ and $B$ do not intersect transversely at $x$;
then we have,
\[
T_{A,x}+T_{B,x}\neq T_{X,x} \]
\[
\dim (T_{A,x}\cap T_{B,x})
> \dim T_{A,x}+\dim T_{B,x} -\dim X
\geq \dim A +\dim B - \dim X \]
Since $W$ is an open subset of $A \cap B$ containing $x$ we have,
\[
\dim T_{W,x} = \dim (T_{A,x}\cap T_{B,x})
> \dim A+\dim B-\dim X = \dim W \]
This contradicts the smoothness of $W$ and demonstrates the lemma. 
\end{proof}

\subsection{Inequality for Families of Subvarieties with Dominant Proper Intersection Locus}\ \\
\indent
We will now begin to consider families of subvarieties in $X$.
It is natural to allow the objects in the families to degenerate to non-integral subschemes.
In particular we will use the following notion of an algebraic family of subschemes of $X$.

\begin{defn}[Algebraic Family of Subschemes]
Given a variety $X$, we will define an
\textit{algebraic family of subschemes} of $X$ of dimension $d$ over a scheme $W$
to be a closed subscheme $\sA\sub X\times W$
such that $p_2:\sA\>W$ is dominant
and every non-empty fiber $\sA_w$
has dimension $d$ over $k(w)$.\\
\indent
If $p_1:\sA\>X$ is dominant we will say $\sA$ is a \textit{covering family}.
If $\sA$ is a variety we will say the family is \textit{integral}.
\end{defn}   

\begin{note}
If $g:[U]\>W$ is a well defined family of proper algebraic cycles of 
a variety $X$ of characteristic $0$
in the sense of \cite[I.3.10]{Koll}
and $U$ is integral,
then $U\sub X\times W$ is an algebraic family of subschemes
and for every $w\in W$,
\[
\supp g^{[-1]}(w)=\supp U_w \]
\end{note}

We will want to consider pairs of covering algebraic families of subschemes which will
correspond to two subvarieties $A$ and $B$ moving with great freedom in $X$.
We will also ask that these families intersect in a``nice'' way defined as follows.

\begin{defn}[Dominant Proper Intersection Locus]
Let $\sA\sub X\times S,\ \sB\sub X\times T$ 
be integral algebraic families of subschemes over varieties $S$ and $T$.
We will say that $\sA$ and $\sB$ have a \textit{dominant proper intersection locus}
if there exists an open subset $U_{\prp}\sub X$ such that
for all closed points $x\in U_{\prp}$, there exists closed points $s\in S,\ t\in T$ 
with $\sA_s$ and $\sB_t$ intersecting properly and $x\in\sA_s \cap \sB_t $.
\end{defn}

For the rest of the discussion in this subsection,
fix two integral algebraic families of subschemes,
$\sA\sub X\times S,\ \sB\sub X\times T$ 
of dimension $d_A$ and $d_B$
over varieties $S$ and $T$
with a dominant proper intersection locus.
Notice that this last condition implies
that both families are covering families
and that $d_A+d_B\geq\dim X$.\\
\indent
Consider the intersection subscheme,
\[
Y=\sA\times T \cap \sB\times S \sub X\times S\times T \]

\begin{prop}
$Y$ has an irreducible component, $Y'$, for which
$p_1:Y'\> X$ and $p_2\times p_3:Y'\> S\times T$ are dominant and for which,
\[
\dim Y'=\dim \sA +\dim \sB -\dim X \]
\end{prop}
\begin{proof}
Since $\sA$ and $\sB$ are dominant over $X$ we know by \cite[ex. II.3.22]{Hart1}
that there is an open set $U'\sub X$ such that for 
$x\in U'\ \dim \sA_x=\dim \sA-\dim X,\ \dim\sB_x=\dim\sB-\dim X$. These must be pure dimensions.
We always have,
\[
Y_x\cong \sA_x\times \sB_x \]
Thus for $x\in U'$, we have,
\[
\dim Y_x = \dim\sA+\dim\sB-2\dim X \]
By definition $\sA_s$ intersects $\sB_t$ properly exactly when 
$\dim Y_{s,t} = \dim \sA_s\cap\sB_t = d_A+d_B-\dim X$
(i.e. the dimension of the intersection is as small as possible).
This means by \cite[ex. II.3.22]{Hart1} again that this is an open condition in the scheme-theoretic image
$p_2\times p_3(Y)\sub S\times T$.
Let this properness condition correspond to $R\sub p_2\times p_3(Y)$.\\
\indent
Then $p_1^{-1}(U'\cap U_{\prp})\cap (p_2\times p_3)^{-1}(R)\sub Y$ is an open subset
and according to our assumption dominates $U'\cap U_{\prp}$.
Let $Y'\sub Y$ be the closure of an 
irreducible component of this open set which dominates $U'\cap U_{\prp}$
and is of maximal dimension with respect to this criterion.
Then dimensional considerations show that,
\[
\dim Y' = \dim \sA+\dim\sB-\dim X \]
and that $Y'$ must also dominate $S\times T$.
\end{proof}

\begin{prop}\label{find W'}
There exists a non-empty open subset $W'\sub Y'\cap (U\times S\times T)$ such that
for all closed points, $(x,s,t)\in W'$,
$\phz(\sA_s)$ and $\phz(\sB_t)$ are smooth at $\phz(x)$.
\end{prop}
\begin{proof}
Consider the open set $p_1^{-1}(U)\sub\sA$ and the subvariety,
\[
(\phz\vert_U\times \id) (p_1^{-1}(U))\sub V\times S \]
Let $Y(\sA)\sub (\phz\vert_U \times \id) (p_1^{-1}(U))$
be an open subset on which the projection to $S$ is smooth
and let $W(\sA)=(\phz\vert_U\times \id)^{-1}(Y(\sA))$.
Then $W(\sA)\sub\sA$ is open and for any closed point $(x,s)\in W(\sA)$,
$\phz(\sA_s)$ is smooth at $\phz(x)$.
Similarly we may define $W(\sB)\sub\sB$.\\
\indent
Since $\sA$ and $\sB$ are integral, we have,
\begin{align*}
\dim W(\sA)^c &\leq \dim \sA-1\\
\dim W(\sB)^c &\leq \dim \sB-1 \end{align*}
This means that for a general $x\in X$,
\begin{align*}
\dim W(\sA)_x^c&\leq\dim\sA-\dim X -1 \\
\dim W(\sB)_x^c&\leq\dim\sB-\dim X -1 \\
\dim W(\sA)_x^c\times \sB_x \cup \dim W(\sB)_x^c\times \sA_x &\leq \dim \sA+\dim\sB-2\dim X -1
\end{align*}
In this case it follows that $Y'\cap (W(\sA)\times T)\cap (W(\sB)\times S)\sub Y'$ is non-empty.
Let $W'$ be this open subset.
\end{proof}

Now let $K=\overline{K(S\times T)}$ be the algebraic closure of the function field of $S\times T$.
Let $X_K, U_K, (U_{\sm})_K$ be the appropriate varieties under field extension
and $A_K, B_K, Y_K, Y'_K$ be the pullbacks of 
$\sA\times T, \sB\times S, Y, Y'$
under the generic map $\Spec K\> S\times T$.
Similarly let $\phz_K:X_K\ratmap V_K$ be the rational map corresponding to $\phz$ under field extension.

\begin{prop}
There is a non-empty open set $W'_K\sub Y'_K\cap U_K$ such that for all closed points $x\in W'_K$,
$\phz_K(A_K)$ and $\phz_K(B_K)$ are smooth at $\phz_K(x)$.
\end{prop}
\begin{proof}
Let $W'\sub Y'$ be as in proposition \ref{find W'}.
If we follow that proof we find that
$(Y(\sA)\times T)_K\sub \phz_K(A_K)$ defines a non-singular open set.
We have,
\[
W(\sA)_K=\phz_K^{-1}((Y(\sA)\times T)_K) \]
and $W(\sB)_K$ is defined similarly.
Then $W'_K = Y'_K\cap W(\sA)_K\cap W(\sB)_K$ has the desired property.
\end{proof}

\begin{prop}\label{find W_K}
There exists a non-empty smooth open subset, $W_K\sub Y_K\cap (U_{\sm})_K$ such that
$\dim_K W_K=\dim_K A_K+\dim_K B_K -\dim_K X_K$ and for all closed points $x\in W_K$,
$\phz_K(A_K)$ and $\phz_K(B_K)$ are smooth at $\phz_K(x)$.
\end{prop}
\begin{proof}
$Y'_K$ is a component of $Y_K$.
Since $\dim Y'=\dim \sA+\dim\sB-\dim X$, it follows that,
\[
\dim_K Y'_K=d_A+d_B-\dim X =\dim_K A_K+\dim_K B_K -\dim_K X_K \]
Let $Y_{\sm}$ be the smooth locus of $Y_K$.
Then let $W_K=Y_{\sm}\cap W'_K\cap (U_{\sm})_K\sub Y'_K$.
It is not hard to check that $W_K$ is non-empty and satisfies the conditions in the proposition.
\end{proof}

\begin{cor}\label{K ineq}
Let $\phz_K:X_K\ratmap V_K$, $A_K$ and $B_K$ be as in this section, then,
\[
\dim V_K\leq \dim \phz_K(A_K)+\dim \phz_K(B_K) \]
\end{cor}
\begin{proof}
This follows from proposition \ref{find W_K}, proposition \ref{dim V leq}
and lemma \ref{trans cond}.
\end{proof}

\subsection{Inequality for Subvarieties Moving as Cycles in Covering Families with Dominant Proper Intersection Locus}
\ \\
\indent
Finally, we will want to use corollary \ref{K ineq}
to prove a relation between Kodaira-Iitaka dimension on a projective variety
and Kodaira-Iitaka dimension on subvarieties that move as cycles.
In order to this we will have to do some work with well defined families of algebraic cycles.

\begin{prop}\label{aKID usc 2}
Let $g:\sA\>S$ be a projective morphism of varieties over $k$.
Suppose $[\sA]\>S$ is a well defined family of algebraic cycles
over $S$ in the sense of \cite[I.3.10]{Koll}.
Let $\sL$ be an invertible sheaf on $\sA$.
Let $\xi\in S$ be the generic point and $s\in S$ an arbitrary point.
Then there exists a subscheme $A \sub g^{-1}(s)$ proper over $k(s)$ such that 
$[A]=g^{[-1]}(s)$ and,
\[
\ka_a(g^{-1}(\xi),\sL_\xi)\leq\ka_a(A,\sL\vert_A) \]
\end{prop}
\begin{proof}
By generic flatness, there is an open subset, $S_0\sub S$
such that $g$ is flat over $S_0$.
Then by proposition \ref{aKID usc}
the function $\ka_a(g^{-1}(y),\sL_y)$ 
is upper semicontinuous for $y\in S_0$.\\
\indent
Choose $h:T\>S$ where $T$ is the spectrum of a DVR, $h(T_g)\in S_0$ and $h(T_0)=s$.
Then according to the definition of a well defined family of algebraic cycles
there exists a subscheme $\sA_{T}'\sub \sA\times_h T$ with $g':\sA_{T}'\> T$ a flat morphism
such that $g'^{-1}(T_g)=g^{-1}(h(T_g))$ and $[g'^{-1}(T_0)]=g^{[-1]}(s)$.
Since $g'$ is flat there is an upper semicontinuity condition for $\ka_a(g'^{-1}(y),\sL_y)$.
Putting these together gives,
\[
\ka_a(g^{-1}(\xi),\sL_\xi)
\leq \ka_a(g^{-1}(h(T_g)),\sL_{h(T_g)})
\leq \ka_a(g'^{-1}(T_0),\sL\vert_{g'^{-1}(T_0)})
\]
Let $A=g'^{-1}(T_0)$.
\end{proof}

\begin{thm}\label{mainII}
Let $X$ be a projective normal variety over an algebraically closed field $k$ of characteristic $0$.
Let $A, B$ be closed normal subvarieties.
Suppose $[A]$ and $[B]$ each move in a well defined integral family of proper algebraic cycles of $X$ over a variety, such that for a general closed point $x\in X$
there are members of the families, $[A'], [B']$,
such that $A'$ and $B'$ intersect properly and $x\in A'\cap B'$.
Then for any invertible sheaf $\sL$ on $X$ we have,
\[
\ka(X,\sL)\leq \ka(A,\sL\vert_A)+\ka(B,\sL\vert_B) \]
\end{thm}
\begin{proof}
Since the families are integral and contain as members the prime cycles $[A]$ and $[B]$,
they will be given by $g_A:[\sA]\> S$ and $g_B:[\sB]\> T$,
where $S$ and $T$ are varieties and $\sA\sub X\times S$ and $\sB\sub X\times T$
are closed subvarieties.
Then $\sA$ and $\sB$ are integral families of subschemes over varieties and
the intersection criterion in the theorem ensures that 
they have a dominant proper intersection locus.
This will allow us to apply corollary \ref{K ineq} at the appropriate place.\\
\indent
If $\ka(X,\sL)=-\infty$ we are done.
Thus assume $\ka(X,\sL)\geq 0$.
Choose $m>0$ so that if $\phz:X\ratmap \P^N$
is the rational map determined by the complete linear system $H^0(X,\sL^m)$
with $V=\phz(X)$, then $\dim V=\ka(X,\sL)$.\\
\indent
By taking pullbacks we may consider $\sL$ to live on $\sA$ and $\sB$ as well as $X$.
Let $K=\overline{K(S\times T)}$ and $\xi_S\in S,\ \xi_T\in T$ be the generic points.
Notice that $A_K$ and $B_K$ in the last section correspond to
$g_A^{-1}(\xi_S)_K$ and $g_B^{-1}(\xi_T)_K$
and the rational map $\phz_K:X_K\ratmap V_K$ comes from the linear system $H^0(X_K,\sL_K^m)$.
Furthermore, since the base locus of $\sL_K^m$
is equal to $(\bl(\sL^m))_K$
and both $\sA$ and $\sB$ are dominant over $X$,
neither $A_K$ nor $B_K$ is contained in $\bl(\sL_K^m)$.
Similarly neither $A_K$ nor $B_K$ is contained in $\Sing X_K$.
Also, for some $s_0\in S$, $A$ is the unique proper subscheme of $X$ for which
$[A]=g_A^{[-1]}(s_0)$.
A similar statement holds for $B$.\\
\indent
Given all of this, we have the following series of inequalities, justified as indicated.
\begin{align*}
\ka(X,\sL)&=\dim_K V_K&\\
&\leq \dim \phz_K(A_K)+\dim \phz_K(B_K)
&\quad{\rm (corollary\ \ref{K ineq})}\\
&\leq \ka_g(A_K,\sL_K)+\ka_g(B_K,\sL_K)
&\quad{\rm (proposition\ \ref{kod im 2})}\\
&\leq \ka_a(A_K,\sL_K)+\ka_a(B_K,\sL_K)
&\quad{\rm (proposition\ \ref{gKID leq aKID})}\\
&= \ka_a(g_A^{-1}(\xi_S),\sL_{\xi_S})+\ka_a(g_B^{-1}(\xi_T),\sL_{\xi_T})
&\quad{\rm (proposition\ \ref{fe aKID})}\\
&\leq \ka_a(A,\sL\vert_A)+\ka_a(B,\sL\vert_B)
&\quad{\rm (proposition\ \ref{aKID usc 2})}\\
&\leq \ka(A,\sL\vert_A)+\ka(B,\sL\vert_B)
&\quad{\rm (proposition\ \ref{KID=gKID=aKID})}\\
\end{align*}
\end{proof}

We will now look at some examples in order to give context to theorem \ref{mainII}.

\begin{exam}\label{exam1}
This example will show that the inequality in theorem \ref{mainII}
may be strict in some circumstances.\\
\indent
Let $X$ be the Hirzebruch surface $\P_{\P^1}(\sO\oplus\sO(-1))$.
The divisor class group of $X$ is generated by divisors $C_0$ and $F$,
where $F$ is a fiber of the natural surjection $\pi:X\>\P^1$
and $C_0$ is the unique section of $\pi$ in $X$ with $C_0^2=-1$.\\
\indent
The divisor $C_0+2F$ is very ample on $X$.
Let $B$ be a general smooth member of the linear system $\vert C_0+2F\vert$
which can be chosen to pass through a general point,
and let $A\sim F$ be a general fiber.
Then $A$ and $B$ move in natural covering families
with dominant proper intersection locus.
Let $\sL=\sO_X(C_0)$. Then,
\begin{align*}
&\deg \sL\vert_A = C_0 .F = 1  &\implies \quad \ka(A,\sL\vert_A)=1  \\
&\deg \sL\vert_B = C_0 .(C_0+2F)=1 &\implies \quad \ka(B,\sL\vert_B)=1 \end{align*}
\indent
However, since $C_0$ is a fixed divisor,
we have $\ka(X,\sL)=0$.
Thus the inequality in theorem \ref{mainII} need not be an equality.
\end{exam}

\begin{exam}
This example will show that some condition on the proper intersection locus of $A$ and $B$
is necessary for the result of theorem \ref{mainII} to hold.
In particular, it is not enough to require that $A$ and $B$ move in covering families
and that they intersect properly.\\
\indent
Let $V=\P^1\times\P^2\sub\P^5$ as in the Veronese embedding.
Let $C(V)\sub\bA^6$ be the cone over $V$
and let $Y$ be the natural closure of $C(V)$ in $\P^6$.
We may blow up the cone point in $Y$ to produce a variety
$\tilde{Y}$ which is a line bundle over $V$ with surjection $\pi:\tilde{Y}\>V$.
$\tilde{Y}$ will be a toric variety.
Its combinatoric construction is described in  \cite[\S 3]{Reid}.
In the same place it is observed that there are contractions of $\tilde{Y}$
to smooth toric varieties produced by projecting the exceptional divisor, $E\iso\P^1\times\P^2$
onto one of its factors.
Let $f:\tilde{Y}\>X$ be the contraction produced by projecting $E$ onto $\P^2$.\\
\indent
Let $\tilde{A}=\pi^{-1}((P_A, H_A)),\ \tilde{B}=\pi^{-1}((P_B, H_A)),\ \tilde{L}=\pi^{-1}((P_L,\P^2))$
for general points $P_A, P_B, P_L\in\P^1$ and general hyperplanes $H_A,H_B\sub \P^2$.
Then for general choices, $\tilde{A}$, $\tilde{B}$, and $\tilde{L}$ will be pairwise disjoint.
However, if we let $A=f(\tilde{A}),\ B=f(\tilde{B}),\ L=f(\tilde{L})$, $A$ and $B$ will meet at a unique
point in the image of $E$ (and thus meet properly)
and each will meet $L$ at a curve contained in $f(E)$.
It is clear that $A$ and $B$ will each move in the same family of cycles covering $X$.
If we let $\sL=\sO_X(L)$, then this invertible sheaf will correspond to a rational map $X\ratmap\P^1$
and it is clear that $\ka(X,\sL)=1$.\\
\indent
By a computation using toric combinatorical methods it can be shown that
$A$ is isomorphic to the Hirzebruch surface in example \ref{exam1}
and that $\sL\vert_A$ corresponds to the divisor $C_0$.
Since this divisor is fixed, we have $\ka(A,\sL\vert_A)=\ka(B,\sL\vert_B)=0$.
All together we have,
\[
1=\ka(X,\sL)>\ka(A,\sL\vert_A)+\ka(B,\sL\vert_B) =0 \]
\indent
This shows that some hypothesis like the dominant proper intersection locus hypothesis
is necessary in theorem \ref{mainII}.
In this example the difficulty is that $A$ and $B$ meet at a point in the base locus of $\sL$.
Perhaps it is possible to find counterexamples to the inequality in theorem \ref{mainII}
where $A$ and $B$ meet at points in the locus where the rational map $\phz_{\sL}$
is defined but is not smooth,
or at points $x\in X$ for which $\phz_{\sL}(x)$ is not a smooth point
of both $\phz_{\sL}(A)$ and $\phz_{\sL}(B)$.
According to proposition \ref{dim V leq} these are the possibilities
to look for in such counterexamples.
\end{exam} 

\bibliographystyle{plain}
\bibliography{KodBib}

\end{document}